\documentclass[12pt]{article}
\usepackage{mathbbold}
\usepackage{amsmath}
\usepackage{mathrsfs}%[10pt]定义字体的大小 {article}指的是文档类别
\usepackage{bbm}
\usepackage{amsfonts}
\usepackage{latexsym,amssymb,amsmath,amscd,amsfonts,epsfig,graphicx}%需要调用的一些软件包
\usepackage[margin=39.5mm]{geometry}

\newtheorem{thm}{Theorem}[section]%定义定理环境
\newtheorem{lem}[thm]{Lemma}%定义性质
\newtheorem{rmk}{Remark}[section]

\newtheorem{con}[thm]{Construction}
\newtheorem{exa}[thm]{\it Example}
\def\pf{\noindent{\it Proof.\;\;\:}}
\def\qed{\nopagebreak\hfill{\rule{4pt}{7pt}}
\medbreak}
\begin{document}

\pagebreak
\begin{center}
{\bf \large A construction of imprimitive symmetric graphs which are not multicovers of their quotients}
\end{center}
\begin{center}
Bin Jia\\
[6pt]

Center for Combinatorics\\
 LPMC - TJKLC, Nankai University\\
 Tianjin 300071\\
P. R. China\\
jiabinqq@gmail.com\\
[5pt]
\end{center}

\noindent {\bf Abstract.} Let $\Sigma$ be a finite $X$-symmetric
graph of valency $\tilde{b} \geq 2$, and $s \geq 1$ an integer. In
this article we give a sufficient and necessary condition for the
existence of a class of finite imprimitive $(X,
s)$-arc-transitive graphs which have a quotient isomorphic to $\Sigma$ and are not multicovers of that quotient, together with a combinatorial method, called the double-star graph construction, for constructing such graphs. Moreover, for any $X$-symmetric graph $\Gamma$ admitting a nontrivial $X$-invariant partition $\mathcal {B}$ such that $\Gamma$ is not a multicover of $\Gamma_{\mathcal {B}}$, we show that there exists a sequence of $m + 1$ $X$-invariant partitions $$\mathcal {B} = \mathcal {B}_0, \mathcal {B}_1, \cdots, \mathcal {B}_m$$ of $V(\Gamma)$, where $m \geq 1$ is an integer, such that $\mathcal {B}_i$ is a proper refinement of $\mathcal {B}_{i - 1}$, $\Gamma_{\mathcal {B}_i}$ is not a multicover of $\Gamma_{\mathcal {B}_{i-1}}$ and $\Gamma_{\mathcal {B}_i}$ can be reconstructed from $\Gamma_{\mathcal {B}_{i-1}}$ by the double-star graph construction, for $ i = 1, 2, \cdots, m$, and that either $\Gamma \cong \Gamma_{\mathcal {B}_m}$ or $\Gamma$ is a multicover of $\Gamma_{\mathcal {B}_m}$.
\vskip 5pt


\section{Introduction}
Let $\Gamma$ be a finite nonempty graph and $s \geq 1$ an integer. A sequence
of $s + 1$ vertices of $\Gamma$ is called an \textit{$s$-arc} if any
two consecutive terms are adjacent and any three consecutive terms
are distinct. A $1$-arc is also called an arc for short. Denote by $Arc_s(\Gamma)$ the set of $s$-arcs of $\Gamma$.
A finite group $X$ acting on the vertices of $\Gamma$
is said to preserve the structure of $\Gamma$ if two vertices of
$\Gamma$ are adjacent if and only if their images under any element
of $X$ are adjacent. If in addition $X$ is transitive on the set of
vertices and transitive on the set of $s$-arcs of $\Gamma$, then
$\Gamma$ is called \textit{$(X, s)$-arc-transitive}. An $(X,
1)$-arc-transitive graph $\Gamma$ is also called an \textit{$X$-symmetric
graph}.

Let $\Gamma$ be a finite $X$-symmetric graph and $\mathcal {B}$ a partition of $V(\Gamma)$. $\mathcal {B}$ is said to be \textit{$X$-invariant} if $B^x \in \mathcal {B}$ for any $B \in \mathcal {B}$ and $x \in X$, where $B^x := \{\sigma^x| \sigma \in B\}$. $\mathcal {B}$ is called \textit{nontrivial} if there exists some $B \in \mathcal {B}$ such that $1 < |B| < |V(\Gamma)|$.

A finite $X$-symmetric graph $\Gamma$ is \textit{imprimitive} if its vertex set admits
a nontrivial $X$-invariant partition $\mathcal {B}$;
otherwise it is called \textit{primitive}. In the imprimitive case,
the \textit{quotient graph} with respect to $\mathcal {B}$ is defined to
have vertex set $\mathcal {B}$ in which two blocks are adjacent
if and only if there exists at least one edge of $\Gamma$ between them. In this article we always assume that $\Gamma_{\mathcal {B}}$ is nonempty, that is, the valency $b$ of $\Gamma_{\mathcal {B}}$ is a positive integer. In this situation, $\Gamma_{\mathcal {B}}$ is $X$-symmetric and all blocks of $\mathcal {B}$ are independent sets of $\Gamma$ (see \cite{AGT}).

Let $\mathcal {B}$ be an $X$-invariant partition of $V(\Gamma)$. For any vertex $\sigma$ of $\Gamma$, denote by $\Gamma(\sigma)$ the set of neighbors of $\sigma$ in $\Gamma$, by $\Gamma_{\mathcal {B}}(B)$ the neighbors of $B$ in $\Gamma_{\mathcal {B}}$, and by $\Gamma_{\mathcal {B}}(\sigma)$ the set of blocks of $\mathcal {B}$ containing at least one neighbor of $\sigma$ in $\Gamma$. Let $b := val(\Gamma_{\mathcal {B}})$, $r := |\Gamma_{\mathcal {B}}(\sigma)|$. Then $b \geq r \geq 1$. For any block $B$ of $\mathcal {B}$, denote by $v$ the number of vertices in $B$, and by $\Gamma(B)$ the set of vertices of $\Gamma$ having at least one neighbor in $B$. For any block $C$ in $\Gamma_\mathcal {B}(B)$, denote by $\Gamma[B, C]$ the bipartite subgraph of $\Gamma$ induced by $(B\cap \Gamma(C))\cup (C \cap\Gamma(B))$. Since $\Gamma$ is $X$-symmetric, $\Gamma[B, C]$ is
$X_{B \cup C}$-symmetric and is independent of the choice of the adjacent blocks $B, C$ of $\mathcal {B}$ up to isomorphism.
Let $d \geq 1$ be the  valency of $\Gamma[B, C]$ and $k := |B \cap \Gamma(C)|$. Then $v \geq k \geq d \geq 1$. Again as $\Gamma$ is $X$-symmetric, the quintuple \textit{$p(\Gamma, X, \mathcal {B})$} $:= (v, k, r, b, d)$ is independent of the choice of $\sigma \in V(\Gamma)$, $B \in \mathcal {B}$ and $C \in \Gamma_{\mathcal {B}}(B)$. Further assume that $\mathcal {B}$ is nontrivial. Then as $vr = bk$, either $v = k \geq 2$ and $r = b \geq 1$, or $v > k \geq 1$ and $b > r \geq 1$. We say $\Gamma$ is or is not a multicover of $\Gamma_{\mathcal {B}}$ respectively in these cases. The second case happens if and only if the subgraph of $\Gamma$ induced by vertices
from $B$ and $C$ is nonempty and contains at least one isolated vertex.

Over the past few decades, finite primitive symmetric graphs have
been studied extensively, and several classification results have
been achieved based on the O'Nan-Scott theorem \cite{LPS} and the
classification of finite simple groups. In contrast, there are not
so many powerful algebraic tools available for dealing with
imprimitive symmetric graphs. In general, any imprimitive symmetric
graph $\Gamma$ has a primitive or bi-primitive quotient graph and
hence the study of the former can be reduced to that of the latter.
However, a lot of important information about $\Gamma$ is lost
during this reduction. Therefore it is difficult to characterize and
reconstruct an imprimitive graph from its quotient graphs.

In \cite{GP}, Gardiner and Praeger proposed an approach to the study of the imprimitive graph $\Gamma$ via investigating the design with point set
$B$ and block set $\{B\cap \Gamma(D)|D \in \Gamma_{\mathcal
{B}}(B)\}$. In \cite{ACOFSGWTTQ},  Li, Praeger and Zhou showed that
for any finite $X$-symmetric graph $\Gamma$ having a nontrivial $X$-invariant partition $\mathcal {B}$ with $k = v
- 1$, if the quotient graph $\Gamma_{\mathcal {B}}$ is $(X, 2)$-arc-transitive, then
$\Gamma$ can be reconstructed from $\Gamma_{\mathcal {B}}$ by the 3-arc graph construction. These work motivated quite a few interesting results about
the structure and construction of imprimitive symmetric graphs,
especially those which are not multicovers of their quotient graphs
\cite{GPZ, Zhou04, Zhou05, Zhou08b, Zhou2002}.

Let $\Sigma$ be a finite $X$-symmetric graph of valency no less than two. In connection with the methods and results above, this paper
aims to answer the following three questions to some extent:
\begin{itemize}
\item[{\bf (Q1)}] When does there exist a finite $(X, s)$-arc-transitive graph $\Gamma$ that has a quotient graph isomorphic to $\Sigma$ but $\Gamma$ is not a multicover of that quotient?

\item[{\bf (Q2)}] If such a $\Gamma$ exists, what information of $\Gamma$ can be obtained from $\Sigma$?

\item[{\bf (Q3)}] How to construct such a $\Gamma$?
\end{itemize}

To answer these questions, we will analyze a certain kind of `local' structures of
$\Sigma$. This structure, called a star, is a set of $l$-arcs starting from the same first vertex of $\Sigma$, where $l \geq 1$
is an integer. The stars will be used to construct imprimitive symmetric graphs from given symmetric graphs in Construction \ref{con_ImSymGraphs}.

\section{Main result}
Let $\Sigma$ be a finite $X$-symmetric graph of valency $\tilde{b} \geq 1$ and $l \geq 1$ an integer. For an $l$-arc $\bbsigma = (\sigma, \sigma_1, \cdots, \sigma_l)$ of
$\Sigma$ and a positive integer $i \leq l$, let $$\bbsigma(i) :=
(\sigma, \sigma_{1}, \cdots, \sigma_i).$$ For any set $S$ of
$l$-arcs of $\Sigma$, let $$S(i) := \{\bbsigma(i)|\bbsigma \in S\}.$$

Let $S$ be a set of $l$-arcs of $\Sigma$ starting from the same
first vertex $\sigma$. $S$ is said to be an
\textit{$(\tilde{r}, l)$-star} of $\Sigma$, where
$\tilde{r} \leq \tilde{b}$ is a positive integer, if
$|S(1)| = \tilde{r}$ and for every integer $i$ with $2 \leq i \leq
l$ and any $(i - 1)$-arc $(\sigma, \sigma_1, \cdots, \sigma_{i -
1})$ of $S(i - 1)$, there are exactly $\tilde{r} - 1$ different vertices
$\sigma_i^1, \sigma_i^2, \cdots, \sigma_i^{\tilde{r} - 1}$ of
$\Sigma(\sigma_{i - 1})$, such that $(\sigma, \sigma_1, \cdots,
\sigma_{i - 1}, \sigma_i^j)$ is an $i$-arc of $S(i)$, for $j = 1,
2, \cdots, \tilde{r} - 1$. $\sigma$ is called the \textit{center} of $S$ and is denoted as \textit{$C_{\Sigma}(S)$}. A $(3, 2)$-star is given in Figure $1$.

\begin{rmk}
As an $l$-arc may contain cycles, the graph spanned by the edges of an $(r, l)$-star does not always have the tree-like structure. An instance can be seen in Example \ref{exa_K2n}.
\end{rmk}

For any $(\tilde{r}, l)$-star $S$ of $\Sigma$ and a positive integer $s \leq l$, it is not difficult to see that $S(s)$ is an $(\tilde{r}, s)$-star of
$\Sigma$. Denote by $X_{S}$ the set-wise stabilizers of $S$ in
$X$. Then $X_{S}$ also acts on $S(s)$. We say that $S$ is
\textit{$(X_{S}, s)$-rank-transitive} if the action of
$X_{S}$ on $S(s)$ is transitive.

Let $S$ and $T$ be two $(\tilde{r}, l)$-stars of $\Sigma$ with centers $\sigma$ and $\tau$, respectively. Then $(S, T)$ is called an \textit{$(\tilde{r}, l)$-double-star} of $\Sigma$ (see Figure $1$) if $(\sigma, \tau) \in S(1)$, $(\tau, \sigma) \in T(1)$ and either $l = 1$ or the following two conditions hold:
\begin{figure}\label{F:DoubleStar}
  \centering
  \includegraphics[width=\textwidth]{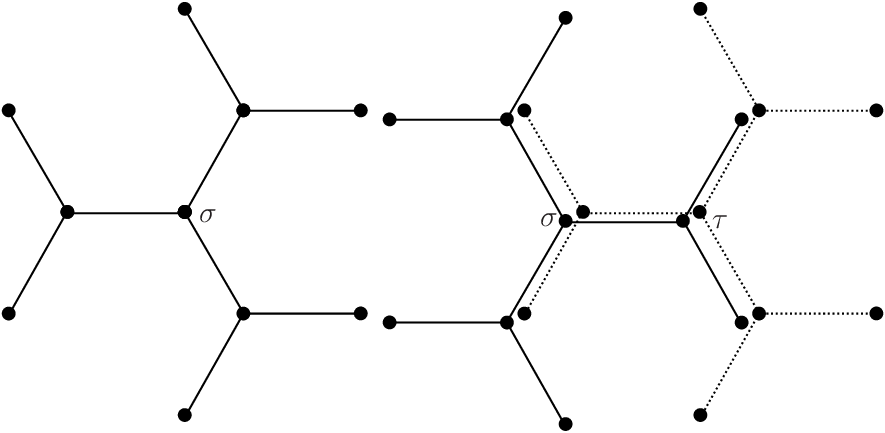}
  \caption{A $(3, 2)$-star and a $(3, 2)$-double-star}
\end{figure}

\begin{itemize}
\item[{\bf (1)}] For any $l$-arc $(\sigma, \tau, \tau_1, \cdots, \tau_{l - 1})$ of $S$ with the first arc $(\sigma, \tau)$, \\$(\tau, \tau_1, \cdots, \tau_{l - 1})$ is an $(l - 1)$-arc of $T(l - 1)$.

\item[{\bf (2)}] For any $l$-arc $(\tau, \sigma, \sigma_1, \cdots, \sigma_{l - 1})$ of $T$ with the first arc $(\tau, \sigma)$, \\$(\sigma, \sigma_1, \cdots, \sigma_{l - 1})$ is an $(l - 1)$-arc of $S(l - 1)$.
\end{itemize}

Let $\Theta$ be a set of $(\tilde{r}, l)$-double-stars of $\Sigma$. $\Theta$ is said to be \textit{self-paired} if $(T, S) \in \Theta$ for any $(S, T) \in \Theta$. Let \textit{$V(\Theta) := \{S, T| (S, T) \in \Theta\}$} and \textit{$E(\Theta) := \{\{S, T\}| (S, T) \in \Theta\}$}. For any positive integer $s \leq l$, we say that $\Theta$ is \textit{$(X, s)$-rank-transitive} if $X$ is transitive on both $\Theta$ and $V(\Theta)$ and $X_{S}$ acts transitively on $S(s)$ for any $S \in V(\Theta)$.

The following example shows that $\Theta$ is not necessarily self-paired even if it is $(X, s)$-rank-transitive for some integer $s \geq 1$.
\begin{exa}
An $(X, s)$-arc-transitive graph $\Sigma$ is said to be $(X, s)$-arc-regular if $|X| = |Arc_s(\Sigma)|$. A cubic $(X, 2)$-arc-regular graph $\Sigma$ is said to be of type $X_2^2$ if any $3$-arc of $\Sigma$ can not be reversed by any $x \in X$ (see \cite{AGT}). In this case, $$\Theta := \{(\{(\sigma, \sigma_1), (\sigma, \tau)\}, \{(\tau, \sigma), (\tau, \tau_1)\})|(\sigma_1, \sigma, \tau, \tau_1) \in Arc_3(\Sigma)\}$$ is a set of $(X, 1)$-rank-transitive $(2, 1)$-double-stars of $\Sigma$ which is not self-paired.
\end{exa}
Now it is ready to state our main result.
\begin{thm}\label{thm_main}
Let $\Sigma$ be a finite $X$-symmetric graph of valency $\tilde{b} \geq 2$ and $\Theta$ a self-paired $(X, s)$-rank-transitive set of $(\tilde{r}, l)$-double-stars of $\Sigma$ such that $X_{S} \cap X_{C_{\Sigma}(T)} = X_{T} \cap X_{C_{\Sigma}(S)}$ for some $(S, T) \in \Theta$, where $l \geq s \geq 1$ and $\tilde{b} > \tilde{r} \geq 1$. Then there exists a finite $(X, s)$-arc-transitive graph $\Gamma$ admitting a nontrivial $X$-invariant partition $\mathcal {B}$ of $V(\Gamma)$ with $d = 1$ and $r = \tilde{r}$ such that $\Gamma_{\mathcal {B}} \cong \Sigma$ but $\Gamma$ is not a multicover of $\Gamma_{\mathcal {B}}$.

Conversely, suppose $\Gamma$ is a finite $(X, s)$-arc-transitive graph admitting a nontrivial $X$-invariant partition $\mathcal {B}$ of $V(\Gamma)$ with $d = 1$ such that $\Gamma$ is not a multicover of $\Gamma_{\mathcal {B}}$. Then there exists a self-paired $(X, s)$-rank-transitive set $\Theta$ of $(r, l)$-double-stars of $\Gamma_{\mathcal {B}}$ such that $X_{S} \cap X_{C_{\Gamma_{\mathcal {B}}}(T)} = X_{T} \cap X_{C_{\Gamma_{\mathcal {B}}}(S)}$ for any $(S, T) \in \Theta$, where $l \geq s \geq 1$.
\end{thm}

\section{A construction of imprimitive graphs}\label{sec_Construction}
In this section, we give the following method for constructing a class of imprimitive  $X$-symmetric graphs from given quotient graphs.
\begin{con}\label{con_ImSymGraphs}
Let $\Sigma$ be a finite $X$-symmetric graph of valency $\tilde{b} \geq 2$ and $\Theta$ a self-paired $(X, s)$-rank-transitive set of $(\tilde{r}, l)$-double-stars of $\Sigma$, where $\tilde{r}$, $s$ and $l$ are positive integers with $\tilde{r} \leq \tilde{b} - 1$ and $s \leq l$. Then the \textit{double-star graph} $\Pi(\Sigma, \Theta)$ of $\Sigma$ with respect to
$\Theta$ is defined to have vertex set $V(\Theta) = \{S, T| (S, T) \in \Theta\}$ and edge set $E(\Theta) = \{\{S, T\}| (S, T) \in \Theta\}$.
\end{con}
\begin{rmk}\label{rmk_ImSymGraphs}
The construction described above is inspired by the $3$-arc graph construction (see \cite{ACOFSGWTTQ}). When $l = 1$, $\tilde{r} = \tilde{b} - 1$ and $\Sigma$ is $(X, 2)$-arc-transitive, one can check that $\Pi$ is isomophic to a $3$-arc graph of $\Sigma$. It can also be seen as a generalization of several existing constructing methods. For example, when $l = 1$, it is the double star graph construction defined in \cite{Jia2008}. And when $\tilde{r} = 1$, which implies $l = 1$, it becomes the arc graph construction (see \cite{Iranmanesh2005}).
\end{rmk}

\begin{lem}
$\Pi := \Pi(\Sigma, \Theta)$ is a finite $X$-symmetric graph with arc set $\Theta$.
\end{lem}

The following Lemma shows that $\Pi$ is imprimitive, $\Sigma$ is isomorphic to a quotient graph of $\Pi$ and $\Pi$ is not a multicover of that quotient.
\begin{lem}\label{lem_ImEq}
For any $\sigma \in V(\Sigma)$, let $V_{\sigma} := \{S \in V(\Theta)| C_{\Sigma}(S) = \sigma \}$. Then $\mathcal {S} := \{V_{\sigma}| \sigma \in V(\Sigma)\}$ is a nontrivial $X$-invariant partition of $V(\Pi)$ with $r = \tilde{r}$ such that $\Pi_{\mathcal {S}} \cong \Sigma$ and $\Pi$ is not a multicover of $\Pi_{\mathcal {S}}$.
\end{lem}
\pf
For any $\sigma \in V(\Sigma)$ and $x \in X$, one can check that $V^x_{\sigma} = V_{\sigma^x}$, so $\mathcal {S}$ is $X$-invariant. Let $(S, T) \in \Theta$ with $C_{\Sigma}(S) = \sigma$ and $C_{\Sigma}(T) = \tau$. Since $\tilde{r} \leq
\tilde{b} - 1$, there exists some $\gamma \in \Sigma(\sigma)$
such that $(\sigma, \gamma) \notin S(1)$. As $\Sigma$ is
$X$-symmetric, there must exist some $x \in X_{\sigma}$ such that
$(\sigma, \tau)^x = (\sigma, \gamma)$. Then $S^x \neq S$
as $(\sigma, \gamma) \in S^x(1)\setminus S(1)$, and $\{S^x,
S\} \subseteq V_{\sigma}$ as $x \in X_{\sigma}$. Note $T \in
V \setminus V_{\sigma}$, so $\mathcal {S}$ is nontrivial. Therefore,
$\Pi$ is a finite $X$-symmetric graph and $\mathcal {S}$ is a nontrivial
$X$-invariant partition of $V(\Pi)$.

Clearly, the mapping:
%\begin{center}
$\delta \mapsto V_{\delta} \in \mathcal {S}$,
for $\delta \in V(\Sigma)$,
%\end{center}
is a bijection between $V(\Sigma)$
and $V(\Pi_{\mathcal {S}})$. On one hand, for any two adjacent
vertices $\sigma_1, \tau_1$ of $\Sigma$, as $\Sigma$ is
$X$-symmetric, there exists some $x \in X$ such that $(\sigma_1,
\tau_1) = (\sigma, \tau)^x$. So $S^x \in V_{\sigma_1}$ and
$T^x \in V_{\tau_1}$ are adjacent in $\Pi$, and hence
$V_{\sigma_1}$ and $V_{\tau_1}$ are adjacent in $\Pi_{\mathcal
{S}}$. On the other hand, for any two adjacent vertices
$V_{\sigma_1}$ and $V_{\tau_1}$ of $\Pi_{\mathcal {S}}$, there
exists some $S_1 \in V_{\sigma_1}$ and $T_1 \in V_{\tau_1}$ such
that $(S_1, T_1) \in \Theta$. Then $(\sigma_1, \tau_1) \in S_1(1)$
is an arc of $\Sigma$. By the analysis above, $\Sigma \cong
\Pi_{\mathcal {S}}$.

For any $(\sigma, \delta) \in S(1)$, as $S$ is $(X_{S},
s)$-rank-transitive, there exists some $x \in X_{S}$ such that
$(\sigma, \delta) = (\sigma, \tau)^x$. Then $(S, T^x) =
(S^x, T^x) \in \Theta$. Therefore, $S$ is adjacent to $T^x \in
V_{\delta}$, and hence $V_{\delta} \in \Pi_{\mathcal {S}}(S)$. Conversely, for any block $V_{\delta}$ of
$\Pi_{\mathcal {S}}(S)$, there exists some $R \in V_{\delta}$ such that
$(S, R)\in \Theta$. Then $\delta$ is the center of $R$. By
the definition of an $(\tilde{r}, l)$-double-star, $(\sigma,
\delta)\in S(1)$. This shows that the mapping:
\begin{center}
$(\sigma, \delta) \mapsto V_{\delta} \in \Pi_{\mathcal {S}}(S)$, for $(\sigma, \delta) \in S(1)$,
\end{center}
is a bijection between $S(1)$ and $\Pi_{\mathcal {S}}(S)$. So $r = |\Pi_{\mathcal {S}}(S)| = |S(1)| = \tilde{r}$ is a positive integer no more than $b - 1$, and thus $\Pi$ is not a multicover of $\Pi_{\mathcal {S}}$.
\qed

The following theorem tells us when $\Pi$ is $(X, s)$-arc-transitive.
\begin{thm}\label{thm_doublestargraphs}
In Construction \ref{con_ImSymGraphs}, if $X_{S} \cap X_{C_{\Sigma}(T)} = X_{T} \cap X_{C_{\Sigma}(S)}$ for any $(S, T) \in \Theta$, then $d = 1$ and $\Pi$ is $(X, s)$-arc-transitive.
\end{thm}
\pf
Let $\sigma = C_{\Sigma}(S)$ and $\tau = C_{\Sigma}(T)$. For any $R \in V_{\tau}$ adjacent to $S$, as $\Pi$ is $X$-symmetric, there exists some $x \in X$ such that $(S^x, T^x) = (S, R)$. As $T$ and $R$ have the same center $\tau$, $x \in X_{S}\cap X_{\tau}= X_{T}\cap X_{\sigma} \leq X_{T}$, thus $R = T^x = T$, and so $d = 1$.

Let $(S_0, S_1, \cdots, S_s)$ be an $s$-arc of $\Pi$. Then $(S_{i-1}, S_i) \in \Theta$ for $i = 1, 2, \cdots, s$; and $S_{i-1} \neq S_{i+1}$ for $i = 1, 2, \cdots, s-1$. Denote by $\sigma_i$ the center of $S_i$ in $\Sigma$, for $i = 0, 1, \cdots, s$. By the definition of an $(\tilde{r}, l)$-double-star, $(\sigma_{i-1}, \sigma_i) \in S_{i-1}(1)$, thus $\sigma_{i - 1}$ is adjacent to $\sigma_i$ in $\Sigma$, for $i = 1, 2, \cdots, s$.

If $\sigma_{i - 1} = \sigma_{i + 1}$, for some $1 \leq i \leq s - 1$, then $S_{i + 1} \in V_{\sigma_{i+1}} = V_{\sigma_{i - 1}}$, and thus $S_i \in V_{\sigma_i}$ is adjacent to both $S_{i-1}$ and $S_{i+1}$ in $\Pi[V_{\sigma_{i - 1}}, V_{\sigma_i}]$, which contradicts $d = 1$. So $\sigma_{i - 1} \neq \sigma_{i + 1}$ for $i = 1, 2, \cdots, s$. Hence $(\sigma_0, \sigma_1, \cdots, \sigma_s)$ is an $s$-arc of $\Sigma$.

Note that $(\sigma_{s - 1}, \sigma_s) \in S_{s - 1}(1)$. Inductively, assume $(\sigma_{s - i}, \sigma_{s - i + 1}, \cdots, \sigma_s) \in S_{s - i}(i)$ for some integer $i$ such that $1 \leq i \leq s - 1$. As $(\sigma_{s - i}, \sigma_{s - i + 1}, \cdots, \sigma_s) \in S_{s - i}(i)$ and $\sigma_{s - i - 1} \neq \sigma_{s - i + 1}$, by the definition of an $(\tilde{r}, l)$-double-star, we have $$(\sigma_{s - i - 1}, \sigma_{s - i}, \sigma_{s - i + 1}, \cdots, \sigma_s) \in S_{s - i - 1}(i + 1).$$ It follows from induction that $(\sigma_0, \sigma_1, \cdots, \sigma_s) \in S_0(s)$.

For any $s$-arc $(L_0, L_1, \cdots, L_s)$ of $\Pi$ starting from $L_0 = S_0$, denote by $\gamma_i$ the center of $L_i$, for $i = 0, 1, \cdots, s$. Then by the analysis above, $\gamma_0 = \sigma_0$ and $(\gamma_0, \gamma_1, \cdots, \gamma_s) \in L_0(s) = S_0(s)$. As $S_0$ is $(X_{S_0}, s)$-rank-transitive, there exists some $x \in X_{S_0}$ such that $(\sigma_0, \sigma_1, \cdots, \sigma_s)^x = (\gamma_0, \gamma_1, \cdots, \gamma_s)$. Note $S_0^x = S_0 = L_0$. For any integer $i$ such that $1 \leq i \leq s$ and $S_{i - 1}^x = L_{i - 1}$, if $S_{i}^x \neq L_{i}$, as $\sigma_i^x = \gamma_i$, then $L_{i - 1}$ is adjacent to both $L_{i}$ and $S_i^x$ in $V_{\gamma_i}$, contradicting $d = 1$. Thus $S_{i}^x = L_{i}$ for $i = 0, 1, \cdots, s$. Again as $\Pi$ is $X$-symmetric, $\Pi$ is $(X, s)$-arc-transitive.
\qed
\begin{rmk}\label{rmk_ImCycles}
In Theorem \ref{thm_doublestargraphs}, further assume that $\tilde{r} = 2$. Then $\Pi$ is a union of cycles.
\end{rmk}
\section{An approach to imprimitive symmetric graphs}
Let $\Gamma$ be a finite $X$-symmetric graph and $\mathcal {B}$ and $\mathcal {B}_1$ two $X$-invariant partitions of $V(\Gamma)$. $\mathcal {B}_1$ is said to be a \textit{refinement} of $\mathcal {B}$, denoted as $\mathcal {B} \geq \mathcal {B}_1$ , if for any $B_1 \in \mathcal {B}_1$, there exists some $B \in \mathcal {B}$ such that $B_1 \subseteq B$. Denote by $\mathcal {B} > \mathcal {B}_1$ when $\mathcal {B}_1$ is a \textit{proper refinement} of $\mathcal {B}$, that is, $\mathcal {B} \geq \mathcal {B}_1$ but $\mathcal {B} \neq \mathcal {B}_1$. When $\mathcal {B} \geq \mathcal {B}_1$, let \textit{$\mathcal {B} / \mathcal {B}_1 := \{\{B_1 \in \mathcal {B}_1| B_1 \subseteq B\}| B \in \mathcal {B}\}$}. Then $\mathcal {B}/ \mathcal {B}_1$ is an $X$-invariant partition of $V(\Gamma_{\mathcal {B}_1})$ such that  $(\Gamma_{\mathcal {B}_1})_{\mathcal {B}/\mathcal {B}_1} \cong \Gamma_{\mathcal {B}}$. Moreover, $\mathcal {B} /\mathcal {B}_1$ is nontrivial if and only if $\mathcal {B} > \mathcal {B}_1$.

Let $\sigma$ be a vertex of $\Gamma$. Denote by $B(\sigma)$ the block of $\mathcal {B}$ containing $\sigma$. For any positive integer $l$, let
\begin{equation*} \label{eq:1}
\begin{split}
Arc_l(\Gamma_{\mathcal {B}}, \sigma)
:= \{(B(\sigma), B(\sigma_1),
\cdots, B(\sigma_l))&\in Arc_l(\Gamma_{\mathcal
{B}})| \\
 (\sigma, \sigma_1, \cdots, \sigma_l) &\in Arc_l(\Gamma)\}.
 \end{split}
 \end{equation*}

Then \textit{$Arc_1(\Gamma_{\mathcal {B}}, \sigma) = \{(B(\sigma), C)| C \in \Gamma_{\mathcal {B}}(\sigma)\}$} is an $(r, 1)$-star of $\Gamma_{\mathcal {B}}$, where $r = |\Gamma_{\mathcal {B}}(\sigma)|$. Let \textit{$\Theta_1 := \{(Arc_1(\Gamma_{\mathcal {B}}, \delta), Arc_1(\Gamma_{\mathcal {B}}, \gamma))|\{\delta, \gamma\} \in E(\Gamma)\}$}. Then $\Theta_1$ is a self-paired $(X, 1)$-rank-transitive set of $(r, 1)$-double-stars of $\Gamma_{\mathcal {B}}$. Clearly, \textit{$\mathcal {B}_1 := \{\{\delta \in V(\Gamma)|Arc_1(\Gamma_{\mathcal {B}}, \delta) = Arc_1(\Gamma_{\mathcal {B}}, \sigma)\}| \sigma \in V(\Gamma)\}$} is an $X$-invariant partition of $V(\Gamma)$ such that $\mathcal {B} \geq \mathcal {B}_1$.

We use $\mathcal {G}$ to denote the set of triples $(\Gamma, X, \mathcal {B})$ such that $\Gamma$ is a finite $X$-symmetric graph, $\mathcal {B}$ is a nontrivial $X$-invariant partition of $V(\Gamma)$ and $\Gamma$ is not a multicover of $\Gamma_{\mathcal {B}}$. Then for any $(\Gamma, X, \mathcal {B}) \in \mathcal {G}$, by \cite{Jia2008}, $\mathcal {B} > \mathcal {B}_1 \geq 1$, $(\Gamma_{\mathcal {B}_1}, X, \mathcal {B}/ \mathcal {B}_1) \in \mathcal {G}$ and  $\Gamma_{\mathcal {B}_1} \cong \Pi(\Gamma_{\mathcal {B}}, \Theta_1)$. Let $p(\Gamma, X, \mathcal {B}_1) = (v_1, k_1, r_1, b_1, d_1)$. Then counting gives that $v_1 | (v, k)$, $r| (r_1, b_1)$, $d_1 | d$, $v_1r_1 = b_1k_1$, $r_1d_1 = rd = val(\Gamma)$, where $(v, k)$ is the greatest common divisor of $v$ and $k$ and $val(\Gamma)$ is the valency of $\Gamma$. Moreover, $p(\Gamma_{\mathcal {B}_1}, X, \mathcal {B}/\mathcal {B}_1) = (v/v_1, k/v_1, r, b, b_1/r)$.

If $\mathcal {B}_1$ is nontrivial and $\Gamma$ is not a multicover of $\Gamma_{\mathcal {B}_1}$, we can apply the above process on $\mathcal {B}_1$ and get a proper refinement $\mathcal {B}_2$ of $\mathcal {B}_1$, and so on. Then we have the following theorem.

\begin{thm}\label{thm_double_star_series}
Let $(\Gamma, X, \mathcal {B}) \in \mathcal {G}$, $\mathcal {B}_0 := \mathcal {B}$ and $$\mathcal {B}_{i+1} := \{\{\delta \in V(\Gamma)| Arc_1(\Gamma_{\mathcal {B}_i}, \delta) = Arc_1(\Gamma_{\mathcal {B}_i}, \sigma)\}|\sigma \in V(\Gamma)\}$$ for $i = 0, 1, \cdots$. Then there exists a unique integer $m \geq 1$ such that $$\mathcal {B} = \mathcal {B}_0
> \mathcal {B}_1 > \cdots
 > \mathcal {B}_{m} = \mathcal {B}_{m+1} = \cdots$$ is a chain of $X$-invariant partitions of $V(\Gamma)$. Let $p(\Gamma, X, \mathcal {B}_i) = (v_i, k_i, r_i, b_i, d_i)$ for $i = 0, 1, \cdots$. Then the following statements hold:
\begin{itemize}
\item[{\bf (1).}] $(\Gamma, X,
\mathcal {B}_i) \in \mathcal {G}$, $v_i \geq 2$ and $1 \leq k_i \leq
v_i - 1$ for $i = 0, 1, \cdots, m - 1$.

\item[{\bf (2).}] $v_{i+1} | (v_i, k_i)$, $r_i|(r_{i+1},
b_{i+1})$, $d_{i+1} |d_i$,
$v_i r_i = b_i k_i$ and $r_i d_i = val(\Gamma)$ for $i = 0, \cdots, m$.

\item[{\bf (3).}] $\Theta_{i + 1} := \{(Arc_1(\Gamma_{\mathcal {B}_i}, \sigma), Arc_1(\Gamma_{\mathcal {B}_i}, \tau))|\{\sigma, \tau\} \in E(\Gamma)\}$ is a self-paired $(X, 1)$-rank-transitive orbit of $(r_i, 1)$-double-stars of $\Gamma_{\mathcal {B}_i}$ such that $\Gamma_{\mathcal {B}_{i+1}} \cong \Pi(\Gamma_{\mathcal {B}_i},
\Theta_{i+1})$ for $i = 0, 1, \ldots, m - 1$.

\item[{\bf (4).}] $(\Gamma_{\mathcal {B}_j}, X, \mathcal {B}_i/\mathcal {B}_j) \in \mathcal
{G}$ with parameters $(v_i/v_j, k_i/v_j, r_i, b_i, b_j/r_i)$ such
that $(\Gamma_{\mathcal {B}_j})_{\mathcal {B}_i/\mathcal {B}_j}
\cong \Gamma_{\mathcal {B}_i}$, for $0 \leq i < j \leq m$.

\item[{\bf (5).}] Either of the following two cases occurs:
\begin{itemize}
\item[{\bf (5.1).}] $v_{m} = k_{m} = 1$, $\mathcal
{B}_{m}$ is trivial and $\Gamma \cong \Pi(\Gamma_{\mathcal {B}_{m - 1}},
\Theta_{m})$.

\item[{\bf (5.2).}] $v_{m} = k_{m} \geq 2$,
$\mathcal {B}_{m}$ is a nontrivial $X$-invariant partition of
$V(\Gamma)$, such that $\Gamma$ is a multicover of $\Gamma_{\mathcal
{B}_{m}}$.
\end{itemize}
\end{itemize}
\end{thm}

\begin{rmk}
Let $(v, k) = p_1^{n_1}\cdots
p_t^{n_t}$, where $n_i$ is a positive integer for $i = 1, \cdots, t$ and $p_1 <
\cdots < p_t$ are primes. Then $1 \leq m \leq \sum_{i = 1}^{t}{n_i} + 1$.
\end{rmk}

Let $B^l(\sigma)
:= \{\delta \in B(\sigma)|\; Arc_l(\Gamma_{\mathcal
{B}}, \delta) = Arc_l(\Gamma_{\mathcal
{B}}, \sigma)\}$. Then for any positive integer $l$, $\mathcal {B}^l := \{B^l(\sigma)|\sigma \in V(\Gamma)\}$ is an $X$-invariant partition of $V(\Gamma)$ and the actions of $X$ on $V_l := \{Arc_l(\Gamma_{\mathcal {B}}, \delta)|\delta \in V(\Gamma)\}$ and that on $\mathcal {B}^l$ are permutationally equivalent under the bijection from $V_l$ to $\mathcal {B}^l$:
\begin{center}
$Arc_l(\Gamma_{\mathcal {B}}, \delta) \mapsto B^l(\delta) \in \mathcal {B}^l$, for $Arc_l(\Gamma_{\mathcal {B}}, \delta)\in V_l$.
\end{center}

The following theorem shows that the subset $Arc_l(\Gamma_{\mathcal {B}}, \sigma)$ of $l$-arcs of $\Gamma_{\mathcal {B}}$ corresponds to a block of $\mathcal {B}_l$ defined in Theorem \ref{thm_double_star_series}.
\begin{thm}\label{thm_PartEqu}
$\mathcal {B}^l = \mathcal {B}_l$ for any integer $l \geq 1$.
\end{thm}
\pf For any vertex $\sigma$ of $\Gamma$, denote by $B_l(\sigma)$ the block of $\mathcal {B}_l$ containing $\sigma$.
Clearly, $\mathcal {B}^1 = \mathcal {B}_1$. Inductively, suppose that $\mathcal {B}^l = \mathcal {B}_l$ holds for some integer $l \geq 1$. For any vertex $\delta$ of $\Gamma$, $Arc_{l+1}(\Gamma_{\mathcal {B}}, \delta) = Arc_{l+1}(\Gamma_{\mathcal {B}}, \sigma)$ if and only if $Arc_{l}(\Gamma_{\mathcal {B}}, \delta) = Arc_l(\Gamma_{\mathcal {B}}, \sigma)$ and for any $\sigma_1 \in \Gamma(\sigma)$ there exists $\delta_1 \in B(\sigma_1) \cap \Gamma(\delta)$ so that $Arc_l(\Gamma_{\mathcal {B}}, \delta_1) = Arc_l(\Gamma_{\mathcal {B}}, \sigma_1)$. In other words, $B^l(\delta) = B^l(\sigma)$ and $\Gamma_{\mathcal {B}^l}(\delta) = \Gamma_{\mathcal {B}^l}(\sigma)$. Since $\mathcal {B}^l = \mathcal {B}_l$, $\sigma$ and $\delta$ are in the same block of $\mathcal {B}^{l+1}$ if and only if they are in the same block of $\mathcal {B}_l$ and they have the same neighbors in $\Gamma_{\mathcal {B}_l}$. By the definition of $\mathcal {B}_{l+1}$, the latter condition is satisfied if and only if $\sigma$ and $\delta$ are in the same block of $\mathcal {B}_{l+1}$. It follows that $\mathcal {B}^{l+1} = \mathcal {B}_{l+1}$. Therefore, the statement of the theorem holds for $l + 1$ and hence by induction for all integers $l \geq 1$.
\qed

For any $\sigma \in V(\Gamma)$, $Arc_l(\Gamma_{\mathcal {B}}, \sigma)$ is a set of $l$-arcs starting from $B(\sigma)$. When $l = 1$, it is an $(r, 1)$-star of $\Gamma_{\mathcal {B}}$. However, things may be different when $l \geq 2$. The following lemma gives a sufficient and necessary condition for $Arc_l(\Gamma_{\mathcal {B}}, \sigma)$ to be an $(r, l)$-star of $\Gamma_{\mathcal {B}}$ for any $l \geq 1$.

\begin{lem}\label{lem_val_Star}
Let $\Gamma$ be a finite $X$-symmetric graph, $\mathcal {B}$ an $X$-invariant partition of $V(\Gamma)$ with $b \geq 1$, $\sigma$ a vertex of $\Gamma$ and $l \geq 1$ an integer. Then $Arc_l(\Gamma_{\mathcal {B}}, \sigma)$ is an $(r, l)$-star of $\Gamma_{\mathcal {B}}$ if and only if $d = d_{l-1}$.
\end{lem}
\pf The statement is trivial for $l = 1$. Suppose it holds for some $l \geq 1$. Then $Arc_{l + 1}(\Gamma_{\mathcal {B}}, \sigma)$ is an $(r, l + 1)$-star of $\Gamma_{\mathcal {B}}$ if and only if for any $\tau \in \Gamma(\sigma)$ and any $\delta \in \Gamma(\sigma) \cap B(\tau)$, both of the following conditions hold:
\begin{itemize}
\item[{ i)}]     $Arc_{l}(\Gamma_{\mathcal {B}}, \delta) = Arc_{l}(\Gamma_{\mathcal {B}}, \tau)$.

\item[{ ii)}]     $Arc_{l}(\Gamma_{\mathcal {B}}, \tau)$ is an $(r, l)$-star of $\Gamma_{\mathcal {B}}$.
\end{itemize}
According to Theorem \ref{thm_PartEqu} and  the definition of $\mathcal {B}^l$, { i)} holds if and only if $\delta \in B_l(\tau)$. By the arbitrariness of $\tau \in \Gamma(\sigma)$ and $\delta \in \Gamma(\sigma) \cap B(\tau)$, this happens if and only if $\Gamma(\sigma) \cap B(\tau) = \Gamma(\sigma) \cap B_l(\tau)$, that is, $d = d_l$.

By the induction hypothesis, { ii)} holds if and only if $d = d_{l - 1}$. Note $d_l | d_{l - 1}$ and $d_{l - 1} | d$, so the lemma holds for $l + 1$ and hence for all $l \geq 1$.
\qed

The following theorem allows us to reconstruct $\Gamma_{\mathcal {B}_l}$ from $\Gamma_{\mathcal {B}}$ and some $(r, l)$-double-stars of $\Gamma_{\mathcal {B}}$.
\begin{thm}\label{thm_Val_DoubleStar}
Let $(\Gamma, X, \mathcal {B}) \in \mathcal {G}$ be such that $d = d_{l - 1}$. Then $$\Theta^l := \{(Arc_l(\Gamma_{\mathcal {B}}, \sigma), Arc_l(\Gamma_{\mathcal {B}}, \tau))|\{\sigma, \tau\} \in E(\Gamma)\}$$ is a set of self-paired $(X, 1)$-rank-transitive $(r, l)$-double-stars of $\Gamma_{\mathcal {B}}$ such that $\Gamma_{\mathcal {B}_l} \cong \Pi(\Gamma_{\mathcal {B}}, \Theta^l)$.
\end{thm}
\pf
For any two adjacent vertices $\sigma$ and $\tau$ of $\Gamma$, by Lemma \ref{lem_val_Star}, both $Arc_l(\Gamma_{\mathcal {B}}, \sigma)$ and $Arc_l(\Gamma_{\mathcal {B}}, \tau)$ are $(r, l)$-stars of $\Gamma_{\mathcal {B}}$ as $d = d_{l - 1}$. One can check that $\Theta^l$ is a self-paired set of $(X, 1)$-rank-transitive $(r, l)$-double-stars of $\Gamma_{\mathcal {B}}$. By Theorem \ref{thm_PartEqu}, the mapping $B_l(\delta) \mapsto Arc_l(\Gamma_{\mathcal {B}}, \delta)$, for $\delta \in V(\Gamma)$, is a bijection between $\Gamma_{\mathcal {B}_l}$ and $\Pi(\Gamma_{\mathcal {B}}, \Theta^l)$. So $\Gamma_{\mathcal {B}_l} \cong \Pi(\Gamma_{\mathcal {B}}, \Theta^l)$.
\qed

The following theorem will be needed in the proof of Theorem \ref{thm_main}.

\begin{thm}\label{thm_d_Equal}
Let $(\Gamma, X, \mathcal {B}) \in \mathcal {G}$, $(\sigma, \tau)$ an arc of $\Gamma$ and $l \geq m$ an integer. Then $X_{Arc_l(\Gamma_{\mathcal {B}}, \sigma)} \cap X_{B(\tau)} = X_{Arc_l(\Gamma_{\mathcal {B}}, \tau)} \cap X_{B(\sigma)}$ if and only if $d = d_m$.
\end{thm}
\pf Let $x$ be any element of $X_{B_{m + 1}(\sigma)}\cap X_{B(\tau)}$. Then $\sigma^x \in B_{m + 1}(\sigma)$ is adjacent to $\tau^x \in B(\tau)$. By the definition of $\mathcal {B}_{m + 1}$, there exists some $\gamma \in B_m(\tau)$ adjacent to $\sigma^x$. So $\tau^x \notin B_m(\tau)$, that is, $x \notin X_{B_m(\tau)} \cap X_{B(\sigma)}$ if and only if $d \neq d_m$. This means that $d = d_m$ if and only if $X_{B_{m + 1}(\sigma)} \cap X_{B(\tau)} \subseteq X_{B_m(\tau)} \cap X_{B(\sigma)}$. By Theorem \ref{thm_double_star_series}, $\mathcal {B}_{m + 1} = \mathcal {B}_m$. Note $X_{B_m(\sigma)} \cap X_{B(\tau)}$ and $X_{B_m(\tau)} \cap X_{B(\sigma)}$ are conjugate under $z \in X$ which reverses $(\sigma, \tau)$, so we have $d = d_m$ if and only if $X_{B_m(\sigma)} \cap X_{B(\tau)} = X_{B_m(\tau)} \cap X_{B(\sigma)}$. By Theorem \ref{thm_double_star_series} and Theorem \ref{thm_PartEqu}, the theorem follows.
\qed

\section{Proof of Theorem \ref{thm_main}}
Let $\Sigma$ be a finite $X$-symmetric graph of valency $\tilde{b} \geq 2$. Suppose there exists a self-paired $(X, s)$-rank-transitive set $\Theta$ of $(\tilde{r}, l)$-double-stars of $\Sigma$ such that $X_{S} \cap X_{C_{\Sigma}(T)} = X_{T} \cap X_{C_{\Sigma}(S)}$ for some $(S, T) \in \Theta$, where $l \geq s \geq 1$ and $\tilde{b} > \tilde{r} \geq 1$. Let $\Gamma := \Pi(\Sigma, \Theta)$, and $\mathcal {B} := \mathcal {S}$, where $\mathcal {S}$ is defined in Lemma \ref{lem_ImEq}. Then by Lemma \ref{lem_ImEq} and Theorem \ref{thm_doublestargraphs}, the first part of Theorem \ref{thm_main} holds.

Now we show the second part. Let $(\Gamma, X, \mathcal {B}) \in \mathcal {G}$ with $d = 1$ such that $\Gamma$ is $(X, s)$-arc-transitive. Let $l = \max{\{m, s\}}$. Then by Theorem \ref{thm_double_star_series}, $d_{l - 1} | d$ and so $d_{l - 1} = d = 1$. And by Theorem \ref{thm_Val_DoubleStar}, $\Theta := \Theta^l$ is a set of self-paired $(X, 1)$-rank-transitive $(r, l)$-double-stars of $\Gamma_{\mathcal {B}}$.

For any arc $(\sigma, \tau)$ of $\Gamma$, let $S := Arc_l(\Gamma_{\mathcal {B}}, \sigma)$ and $T := Arc_l(\Gamma_{\mathcal {B}}, \tau)$. Then $(S, T) \in \Theta$, $C_{\Gamma_{\mathcal {B}}}(S) = B(\sigma)$ and $C_{\Gamma_{\mathcal {B}}}(T) = B(\tau)$. By Theorem \ref{thm_d_Equal}, we have $X_{S} \cap X_{C_{\Gamma_{\mathcal {B}}}(T)} = X_{T} \cap X_{C_{\Gamma_{\mathcal {B}}}(S)}$.

Let $(B(\sigma), B(\sigma_1), \cdots, B(\sigma_s))$ and $(B(\sigma), B(\delta_1), \cdots, B(\delta_s))$ be two $s$-arcs of $S(s)$, where $(\sigma, \sigma_1, \cdots, \sigma_s)$ and $(\sigma, \delta_1, \cdots, \delta_s)$ are two $s$-arcs of $\Gamma$. As $\Gamma$ is $(X, s)$-arc-transitive, there exists some $x \in X_{\sigma} \leq X_{B_m(\sigma)}$ such that $\delta_i^x = \sigma_i$ for $i = 1, \cdots, s$, and hence $(B(\sigma), B(\sigma_1), \cdots, B(\sigma_s)) = (B(\sigma), B(\delta_1), \cdots, B(\delta_s))^x$. By Theorem \ref{thm_double_star_series}, $B_m(\sigma) = B_l(\sigma)$ as $m \leq l$. By the analysis in and above Theorem \ref{thm_PartEqu}, $X_{B_l(\sigma)} = X_{S}$. So $x \in X_{S}$ and hence $\Theta$ is $(X, s)$-rank-transitive.
\section{Examples}
In this section, we give several examples to show the construction of an imprimitive $(X, s)$-arc-transitive graph $\Pi$ from stars and double-stars of an $X$-symmetric graph $\Sigma$ such that $\Pi$ has a quotient graph isomorphic to $\Sigma$ but is not a multicover of that quotient, where $s$ is a positive integer.

For any integer $n \geq 1$, let $[n] := \{1,
2, \cdots, n\}$. Let $V_1 = \{\sigma_i| i \in [n]\}$ and $V_2 = \{\tau_i| i \in [n]\}$ be two disjoint sets of vertices.

The first example shows that $a$ finite $(X, 3)$-arc-transitive graph may admit an $(X, 2)$-arc-transitive quotient that is not $(X, 3)$-arc-transitive. Note the stars defined below do not have the tree-like structure as they contain cycles.
\begin{exa}\label{exa_K2n}
For any integer $n \geq 3$, denote by $K_{2n}$ the complete graph of vertex set $V_1 \cup V_2$. Let $X = Aut(K_{2n}) \cong Sym([2n])$. Then $K_{2n}$ is $(X, 2)$-arc-transitive but not $(X, 3)$-arc-transitive. Let $$S = \{(\sigma_1, \tau_i, \sigma_j, \tau_k)|i \in [n], j \in [n]\setminus \{1\}, k \in [n] \setminus \{i\}\},$$
$$T = \{(\tau_1, \sigma_i, \tau_j, \sigma_k)|i \in [n], j \in [n]\setminus \{1\}, k \in [n] \setminus \{i\}\}$$ and $\Theta = \{(S, T)^x | x \in X\}$. Then $\Theta$ is a self-paired set of $(X, 3)$-rank-transitive $(n, 3)$-double-stars of $K_{2n}$.  Let $\Pi = \Pi(K_{2n}, \Theta)$, $B = \{S^x | x \in X_{\sigma_1}\}$ and $\mathcal {B} = \{B^x | x \in X\}$. Then $(\Pi, X, \mathcal {B})\in \mathcal {G}$ with $m = 1$ and $(v, k, r, b, d) = ({{2n-1}\choose{n}}, {{2n-2}\choose{n-1}}, n, 2n - 1, 1)$. Moreover $\Pi \cong \frac{1}{2}{{2n}\choose{n}}K_{n,n}$ is
$(X, 3)$-arc-transitive of order $2n{{2n-1}\choose{n}}$.
\end{exa}

For any integer $n \geq 3$ and $1 \leq m \leq n$, let $[n]^{(m)}$ be the set of  $m$-subsets
of $[n]$. Then the odd graph $O_n$ (see \cite{AGT}) is defined to have vertex set $[2n - 1]^{(n - 1)}$ such that two vertices are
adjacent if and only if the corresponding subsets are disjoint.

For any connected $(X, 3)$-arc-transitive graph $\Sigma$ of valency no less than three, it has been given in Corollary 4.8 of \cite{FSGTTQ} that $\Sigma$ is a quotient of at least one connected $X$-symmetric but not $(X, 2)$-arc-transitive graph. Now the following example will show that $\Sigma$ may also be a quotient of an $(X, 2)$-arc-transitive graph.

\begin{exa}\label{ex_O4}
Let $X = A_7 < Aut(O_4)$. Then $O_4$ is $(X, 3)$-arc-transitive (see \cite{Li2001}). For the sake of convenience, we use $\{ijk\}$ to denote the subset $\{i, j, k\}$ of $\{1, 2, \cdots, 7\}$. Set
$$S := \left\{\begin{array}{llll}
(\{123\}, \{456\}, \{127\})\\
(\{123\}, \{456\}, \{137\})\\
(\{123\}, \{457\}, \{126\})\\
(\{123\}, \{457\}, \{136\})\\
(\{123\}, \{567\}, \{124\})\\
(\{123\}, \{567\}, \{134\})\end{array}\right\}; T :=
\left\{\begin{array}{llll}
(\{456\}, \{123\}, \{457\})\\
(\{456\}, \{123\}, \{567\})\\
(\{456\}, \{127\}, \{345\})\\
(\{456\}, \{127\}, \{356\})\\
(\{456\}, \{137\}, \{245\})\\
(\{456\}, \{137\}, \{256\})\end{array}\right\}.$$
Let $\Theta = \{(S, T)^x|x\in X\}$. Then $\Theta$ is a self-paired $(X, 2)$-rank-transitive orbit of $(3, 2)$-double-stars of $O_4$. Let $\Pi = \Pi(O_4, \Theta)$, $B = \{S^x | x \in X_{\{123\}}\}$ and $\mathcal {B} = \{B^x | x \in X\}$. Then $(\Pi, X, \mathcal {B}) \in \mathcal {G}$ with $m = 2$, $(v, k, r, b, d) = (12, 9, 3, 4, 1)$ and $(v_1, k_1, r_1, b_1, d_1) = (3, 1, 3, 9, 1)$. Moreover, $\Pi$ is $(X, 2)$-arc-transitive of order $420$.
\end{exa}

For any finite group $X$ and any $x\in X$, denote by $\bar{x}$ the right regular representation of $x$. For any subgroup $G$ of $X$, let $\bar{G} = \{\bar{g}| g \in G\} \leq \bar{X}$.
\begin{exa}\label{ex_PSL2p}
For any prime $p \geq 11$ such that $p \equiv \pm 1 \;(mod\; 8)$, let $X = PSL(2, p)$, and $H$ the subgroup of $X$ isomorphic to $Sym([4])$. By \cite{OPPGWSSATOG}, there exists an involution $z \in X \setminus H$ such that $N_X(P) = P \times <z>$, where $P = H \cap H^z \cong Sym([3])$. Let $\Sigma = Cos(X, H, HzH)$. Then $\Sigma$ is a finite $(\bar{X}, 2)$-arc-transitive graph of valency $4$. Let $[H:P] = \{Ph_0, Ph_1, Ph_2, Ph_3\}$, where $h_0 = 1$. Then $Q_i := P\cap P^{h_i}$, for $i = 1, 2, 3$ are the only three subgroups of $P$ isomorphic to $Z_2$. Since $z \in N_X(P)$, $z$ fixes $\{Q_1, Q_2, Q_3\}$ set-wise by conjugation. Note $z$ is an involution, so $z$ belongs to at least one of $N_X(Q_i)$, for $i = 1, 2, 3$. Without loss of generality, suppose $Q_3^z = Q_3$. Let $$S = \{(H, Hzh_i, Hzh_jzh_i)|i \in \{0, 1, 2\}, j \in \{1, 2\}\}.$$ Then $\Theta := \{(S, S^{\bar{z}})^{\bar{x}}|\bar{x} \in \bar{X}\}$ is a self-paired $(\bar{X}, 2)$-rank-transitive orbit of $(3, 2)$-double-stars of $\Sigma$. Let $\Pi = \Pi(\Sigma, \Theta)$, $B = \{S^{\bar{h}}|h \in H\}$ and $\mathcal {B} = \{B^{\bar{x}}| x \in X\}$. Then $(\Pi, \bar{X}, \mathcal {B}) \in \mathcal {G}$ with $m = 1$ and $(v, k, r, b, d) = (4, 3, 3, 4, 1)$. Moreover, $\Pi$ is $(\bar{X}, 2)$-arc-transitive.
\end{exa}

\end{document}